\documentclass[reqno,centertags,a4paper,10pt]{amsart}
\usepackage{amssymb,enumerate}
\usepackage{amsthm}
\usepackage{eucal}
\usepackage{fullpage}

\usepackage{pslatex}			

\usepackage{lastpage}			
\usepackage{indentfirst}		
\usepackage{color}				
\usepackage{graphicx}			
\usepackage{microtype} 			
\usepackage{amssymb, amsthm}
\usepackage{amsmath}

\usepackage{color}				

\usepackage{hyperref}     

\usepackage{latexsym}
\usepackage[latin1]{inputenc}
\usepackage{multirow}

\usepackage[all]{xy}

\usepackage{pgf,tikz}

\usetikzlibrary{arrows}

\newcommand{\R}{\textnormal{I}\!\textnormal{R}}
\newcommand{\N}{\textnormal{I}\!\textnormal{N}}

\newcommand{\cjap}[1]{\langle {#1}\rangle}

\numberwithin{equation}{section}

\newtheorem{theorem}{Theorem}[section]

\newtheorem{remark}[theorem]{Remark}
\newtheorem{lemma}[theorem]{Lemma}

\begin{document}
	\title[ill-posedness]
	{A note on $C^2$ Ill-posedness results for the Zakharov system in arbitrary dimension}
	
	\author[L. Domingues, R. Santos]
	{Leandro Domingues\ \ \ \&\ \ \ Raphael Santos } 
	
	\address{Leandro Domingues \newline\indent 
	Departamento de Matemática Aplicada, CEUNES/UFES. Rodovia BR 101 Norte, Km 60, Bairro Litorâneo,\newline\indent CEP 29932-540, São Mateus, ES, Brazil
		}
	\email{leandro.ceunes@gmail.com, leandro.domingues@ufes.br}
	
	\address{Raphael Santos \newline\indent   
		UFRJ-Campus Prof. Alo\'isio Teixeira, 
		Granja dos Cavaleiros - Maca\'e,            
		CEP 27930-560, RJ - Brazil }                  
	\email{raphaelsantos@macae.ufrj.br}

	
	\begin{abstract}
		This work is concerned with the Cauchy problem for a Zakharov system with initial data in Sobolev spaces $H^k(\R ^d)\!\times\!H^l(\R ^d)\!\times\!H^{l\!-\!1}\!(\R ^d)$.
		%
		%
		We recall the well-posedness and ill-posedness results known to date and establish new ill-posedness results. We prove $C^2$ ill-posedness for some new indices $(k,l)\in\R^2$. Moreover, our results are valid in arbitrary dimension. 
		We believe that our detailed proofs are built on a methodical approach and can be adapted to obtain similar results for other systems and equations.
	\end{abstract}
	
	\maketitle
	
	\allowdisplaybreaks
	
	\section{Introduction}

	This work is concerned with the Cauchy problem for the following Zakharov system
	\begin{equation}\label{Z}\tag{Z}
		\left\{\begin{gathered}
			i\partial_t u\ +\ \Delta u\ =\ nu\ ,
			\quad \hfill 
			u:\R\!\times\!\R^d\to\mathbb C\,,\\ 
			\partial_t^2n\ -\ \Delta n\ =\ \Delta |u|^2, 
			\quad\quad\quad\quad \hfill
			\ n:\R\!\times\!\R^d\to\R,\\ 
			(u,n,\partial_tn)|_{t=0}\,\in\,H^{k,l}, \hfill\  
		\end{gathered}\right.
	\end{equation}
	where $H^{k,l}$ is a short notation for the Sobolev space $H^k\!(\R ^d;\mathbb C)\!\!\times\!H^l\!(\R ^d;\R)\!\!\times\!H^{l\!-\!1}\!(\R ^d;\R)$, $(k,l)\in\R^2$ and $\Delta$ is the laplacian operator for the spatial variable. 

	V. E. Zakharov introduced the system \eqref{Z} in \cite{zakharov} to describe the long wave Langmuir turbulence in a plasma. The function $u$ represents the slowly varying envelope of the rapidly oscillating electric field and the function $n$ denotes the deviation of the ion density from its mean value.
%

In this note we prove that, for any dimension $d$, the system \eqref{Z} is $C^2$ ill-posed in $H^{k,l}$, 
for the indices $(k,l)$ displayed in Figure \ref{Region-StnotC2} and Figure \ref{Region-SnotC2} 
(see Theorem \ref{main02} and Theorem \ref{main01} for the precise statements). 
The first \emph{$C^2$ ill-posedness} result was proved by Tzvetkov in \cite{Tzvet} for the KdV equation, 
improving the previous \emph{$C^3$ ill-posedness} result of Bourgain found in \cite{Bo2}. 
We essentially follow the same ideas of \cite{Tzvet}, but our proofs are structured as in \cite{leandro}. 
Two slightly different senses of $C^2$ ill-posedness are considered in our results (see also Remark \ref{rmksenseill}). 

\begin{center}
\definecolor{cinza 1}{rgb}{0.875,0.875,0.875}
\definecolor{cinza 2}{rgb}{0.5,0.5,0.5}
\definecolor{cinza 3}{rgb}{0.125,0.125,0.125}

\begin{figure}[htp]
\begin{minipage}[b]{0.45\linewidth}
\centering 
\begin{tikzpicture}[scale=1][line cap=round,line join=round,>=triangle 45,x=1.0cm,y=1.0cm]
\draw [color=cinza 2,dotted, xstep=0.5cm,ystep=0.5cm] (-0.5,-1.0) grid (3.0,2.5);
\draw[->,color=black] (-0.5,0) -- (3.0,0);
\foreach \x in {1,2,3,4,5}
\draw[shift={(\x/2,0)},color=black] (0pt,1pt) -- (0pt,-1pt) node[above] {\tiny $\x$};
\draw[color=black] (2.7,0) node [anchor=south west] {\footnotesize $k$};
\draw[->,color=black] (0,-1) -- (0,2.5);
\foreach \y in {-2,-1,-0.5,1,2,3,4}
\draw[shift={(0.1,\y/2)},color=black] (1pt,0pt) -- (-1pt,0pt) node[left] {\tiny $\y$};
\draw[color=black] (0,2.35) node [anchor=west] {\footnotesize $l$};
\clip(-0.5,-1.0) rectangle (3.0,2.5);
\fill[fill=purple,fill opacity=0.4] (0,-0.25) -- (3.0,-0.25) -- (3.0,-1.0) -- (-0.5,-1.0) -- (-0.5,2.5) -- (1.375,2.5) -- cycle;
\draw [line width=1.1][dash pattern=on 1.5pt off 1.5pt](0,-0.25)-- (3.0,-0.25);
\draw [line width=1.1][dash pattern=on 1.5pt off 1.5pt](0,-0.25)-- (1.375,2.5);
\draw (0.8,1.8) node [right]{$\rotatebox{65}{\tiny $l=2k-1/2$}$};
\draw (1.95,-0.36) node [anchor=south west] {\tiny $l=-1/2$};
\end{tikzpicture}
\caption{$S^{\,t}\!$ is not $C^2$. Theorem \ref{main02}}\label{Region-StnotC2}
\end{minipage}
\begin{minipage}[b]{0.45\linewidth}
\centering 
\begin{tikzpicture}[scale=1][line cap=round,line join=round,>=triangle 45,x=1.0cm,y=1.0cm]
\draw [color=cinza 2,dotted, xstep=0.5cm,ystep=0.5cm] (-0.5,-1.0) grid (3.0,2.5);
\draw[->,color=black] (-0.5,0) -- (3.0,0);
\foreach \x in {1,2,3,4,5}
\draw[shift={(\x/2,0)},color=black] (0pt,1pt) -- (0pt,-1pt) node[below] {\tiny $\x$};
\draw[color=black] (2.7,0) node [anchor=south west] {\footnotesize $k$};
\draw[->,color=black] (0,-1) -- (0,2.5);
\foreach \y in {-2,-1,1,2,3,4}
\draw[shift={(0,\y/2)},color=black] (1pt,0pt) -- (-1pt,0pt) node[left] {\tiny $\y$};
\draw[color=black] (0,2.35) node [anchor=west] {\footnotesize $l$};
\clip(-0.5,-1.0) rectangle (3.0,2.5);
\fill[fill=red,fill opacity=0.35] (-0.5,0) -- (2,2.5) -- (-0.5,2.5) -- cycle;
\fill[fill=red,fill opacity=0.35] (0,-1) -- (3.0,2.0) -- (3,-1) -- cycle;
\draw [line width=1.1][dash pattern=on 1.5pt off 1.5pt](0,-1)-- (3.0,2.0);
\draw [line width=1.1][dash pattern=on 1.5pt off 1.5pt](-0.5,0)-- (2,2.5);
\draw (1.35,2.08) node [right]{$\rotatebox{45}{\tiny $l=k+1$}$};
\draw (2.05,1.72) node [right]{$\rotatebox{45}{\tiny $l=k-2$}$};
\end{tikzpicture}
\caption{$S$ is not $C^2$. Theorem \ref{main01}}\label{Region-SnotC2}
\end{minipage}
\end{figure}
\end{center}
%
%

Ginibre, Tsutsumi and Velo introduced in \cite{GTV97} a heuristic critical regularity for the system \eqref{Z}, 
which is given by $(k,l)\!=\!(\,d/2-3/2\,,\,d/2-2)$. 
In particular, our result in Theorem \ref{main02} with $d=3$ (physical dimension) shows that 
the critical regularity $(0,-1/2)$ is the endpoint for achieving well-posedness by fixed point procedure. 
We point out that local well-posedness at critical regularity is an open problem for $d\ge3$ 
(see Table \ref{bestresultstodate} bellow). 

	The system \eqref{Z} has been studied in several works. 
	Bourgain and Colliander proved in \cite{boco} local well-posedness in the energy norm for $d=2,3$. 
	They construct local solutions applying the contraction principle in 
	$X^{s,b}$ 
	spaces introduced in \cite{Bo}. 
	Local well-posedness in arbitrary dimension under weaker regularity assumptions was obtained in \cite{GTV97} by Ginibre, Tsutsumi and Velo. We recall the last result in the next theorem.  
	\begin{theorem}\label{gtv1}
	(Ginibre, Tsutsumi and Velo \cite{GTV97})
		Let $d\geq 1$. The system \eqref{Z} is locally well-posed, 
		provided 
		\begin{equation}\label{regionlwp}
			\begin{array}{lll}
				\!\!-1/2< k-l\leq 1,\quad &2k\geq l+1/2\geq 0,&\textrm{for} \ d=1\\
				l\leq k\leq l+1, & &\textrm{for all}\  d\geq 2\\
				l\geq 0, &\quad 2k-(l+1)\geq 0, &\textrm{for}\ d=2,3\\
				l> d/2-2,& 2k-(l+1)>d/2-2,\quad &\textrm{for all}\ d\geq 4.
				
				\end{array}
						\end{equation}
								
	\end{theorem}


\begin{center}
\definecolor{cinza 1}{rgb}{0.875,0.875,0.875}
\definecolor{cinza 2}{rgb}{0.5,0.5,0.5}
\definecolor{cinza 3}{rgb}{0.125,0.125,0.125}

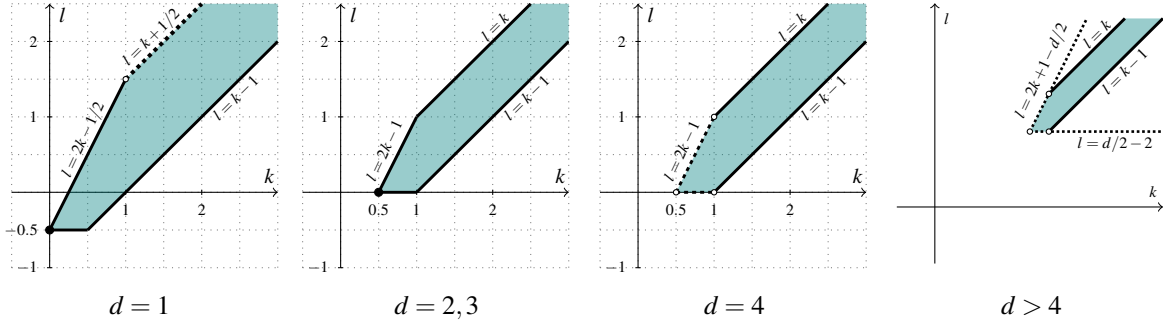
\begin{figure}[htp]
\begin{minipage}[b]{0.24\linewidth}
\begin{tikzpicture}[scale=1][line cap=round,line join=round,>=triangle 45,x=1.0cm,y=1.0cm]
\draw [color=cinza 2,dotted, xstep=0.5cm,ystep=0.5cm] (-0.5,-1.0) grid (3.0,2.5);
\draw[->,color=black] (-0.5,0) -- (3.0,0);
\foreach \x in {1,2}
\draw[shift={(\x,0)},color=black] (0pt,1pt) -- (0pt,-1pt) node[below] {\tiny $\x$};
\draw[color=black] (2.7,0) node [anchor=south west] {\footnotesize $k$};
\draw[->,color=black] (0,-1) -- (0,2.5);
\foreach \y in {-1,-0.5,1,2}
\draw[shift={(0,\y)},color=black] (1pt,0pt) -- (-1pt,0pt) node[left] {\tiny $\y$};
\draw[color=black] (0,2.35) node [anchor=west] {\footnotesize $l$};
\fill[fill=teal,fill opacity=0.4] (0,-0.5) -- (0.5,-0.5) -- (3.0,2.0) -- (3.0,2.5) -- (2,2.5) -- (1,1.5) -- cycle;
\draw [line width=1.1](0.5,-0.5)-- (3.0,2.0);
\draw [line width=1.5][dash pattern=on 1.5pt off 1.5pt] (1,1.5)-- (2,2.5);
\draw [line width=1.1](0,-0.5)-- (0.5,-0.5);
\draw [line width=1.1](0,-0.5)-- (1,1.5);
\draw [fill=black] (0,-0.5) circle (1.5pt);
\draw [fill=white] (1,1.5) circle (1pt);
\draw (2,1.25) node [right]{$\rotatebox{45}{\tiny $l=k-1$}$};
\draw (0.8,2.05) node [right]{$\rotatebox{45}{\tiny $l=k+1/2$}$};
\draw (0,0.7) node [right]{$\rotatebox{65}{\tiny $l=2k-1/2$}$};
\end{tikzpicture}
\centering{$d=1$}
\end{minipage}
\begin{minipage}[b]{0.24\linewidth}
\begin{tikzpicture}[scale=1][line cap=round,line join=round,>=triangle 45,x=1.0cm,y=1.0cm]
\draw [color=cinza 2,dotted, xstep=0.5cm,ystep=0.5cm] (-0.5,-1.0) grid (3.0,2.5);
\draw[->,color=black] (-0.5,0) -- (3.0,0);
\foreach \x in {0.5,1,2}
\draw[shift={(\x,0)},color=black] (0pt,1pt) -- (0pt,-1pt) node[below] {\tiny $\x$};
\draw[color=black] (2.7,0) node [anchor=south west] {\footnotesize $k$};
\draw[->,color=black] (0,-1) -- (0,2.5);
\foreach \y in {-1,1,2}
\draw[shift={(0,\y)},color=black] (1pt,0pt) -- (-1pt,0pt) node[left] {\tiny $\y$};
\draw[color=black] (0,2.35) node [anchor=west] {\footnotesize $l$};
\clip(-0.5,-1.0) rectangle (3.0,2.5);
\fill[fill=teal,fill opacity=0.4] (0.5,0) -- (1,0) -- (3,2) -- (3,2.5) -- (2.5,2.5) -- (1,1) -- cycle;
\draw [line width=1.1](1,0)-- (3.0,2.0);
\draw [line width=1.1](1,1)-- (2.5,2.5);
\draw [line width=1.1](0.5,0)-- (1,1);
\draw [line width=1.1](0.5,0)-- (1,0);
\draw [fill=black] (0.5,0) circle (1.5pt);
\draw (1.7,2.2) node [right]{$\rotatebox{45}{\tiny $l=k$}$};
\draw (2,1.25) node [right]{$\rotatebox{45}{\tiny $l=k-1$}$};
\draw (0.2,0.55) node [right]{$\rotatebox{65}{\tiny $l=2k-1$}$};
\end{tikzpicture}
\centering{$d=2,3$}
\end{minipage}
\begin{minipage}[b]{0.24\linewidth}
\begin{tikzpicture}[scale=1][line cap=round,line join=round,>=triangle 45,x=1.0cm,y=1.0cm]
\draw [color=cinza 2,dotted, xstep=0.5cm,ystep=0.5cm] (-0.5,-1.0) grid (3.0,2.5);
\draw[->,color=black] (-0.5,0) -- (3.0,0);
\foreach \x in {0.5,1,2}
\draw[shift={(\x,0)},color=black] (0pt,1pt) -- (0pt,-1pt) node[below] {\tiny $\x$};
\draw[color=black] (2.7,0) node [anchor=south west] {\footnotesize $k$};
\draw[->,color=black] (0,-1) -- (0,2.5);
\foreach \y in {-1,1,2}
\draw[shift={(0,\y)},color=black] (1pt,0pt) -- (-1pt,0pt) node[left] {\tiny $\y$};
\draw[color=black] (0,2.35) node [anchor=west] {\footnotesize $l$};
\clip(-0.5,-1.0) rectangle (3.0,2.5);
\fill[fill=teal,fill opacity=0.4] (0.5,0) -- (1,0) -- (3,2) -- (3,2.5) -- (2.5,2.5) -- (1,1) -- cycle;
\draw [line width=1.1](1,0)-- (3.0,2.0);
\draw [line width=1.1](1,1)-- (2.5,2.5);
\draw [line width=1.1][dash pattern=on 1.5pt off 1.5pt](0.5,0)-- (1,1);
\draw [line width=1.1][dash pattern=on 1.5pt off 1.5pt](0.5,0)-- (1,0);
\draw [fill=white] (0.5,0) circle (1pt);
\draw [fill=white] (1,0) circle (1pt);
\draw [fill=white] (1,1) circle (1pt);
\draw (1.7,2.2) node [right]{$\rotatebox{45}{\tiny $l=k$}$};
\draw (2,1.25) node [right]{$\rotatebox{45}{\tiny $l=k-1$}$};
\draw (0.2,0.55) node [right]{$\rotatebox{65}{\tiny $l=2k-1$}$};
\end{tikzpicture}
\centering{$d=4$}
\end{minipage}
\begin{minipage}[b]{0.24\linewidth}
\begin{tikzpicture}[scale=1][line cap=round,line join=round,>=triangle 45,x=1.0cm,y=1.0cm]
\draw[color=black] (2.7,0) node [anchor=south west] {\tiny $k$};
\draw[->,color=black] (-0.5,0) -- (3.0,0);
\draw[color=black] (0,2.5) node [anchor=west] {\tiny $l$};
\draw[->,color=black] (0,-0.75) -- (0,2.7);
\clip(-0.5,-1.0) rectangle (3.0,2.7);
\fill[fill=teal,fill opacity=0.4] (1.25,1) -- (1.5,1) -- (3,2.5) -- (2.5,2.5) -- (1.5,1.5) -- cycle;
\draw [line width=1.1](1.5,1)-- (3,2.5);
\draw [line width=1.1](1.5,1.5)-- (2.5,2.5);
\draw [line width=1][dash pattern=on 1pt off 1pt](1.25,1)-- (2,2.5);
\draw [line width=1][dash pattern=on 1pt off 1pt](1.25,1)-- (3,1);
\draw (1.75,0.6) node [anchor=south west] {\tiny $l=d/2-2$};
\draw (0.9,1.8) node[right]{$\rotatebox{65}{\tiny$l=2k+\!1\!-d/2$}$};
\draw (2.05,1.8) node [right]{$\rotatebox{45}{\tiny $l=k-1$}$};
\draw (1.8,2.3) node [right]{$\rotatebox{45}{\tiny $l=k$}$};
\draw [fill=white] (1.25,1) circle (1pt);
\draw [fill=white] (1.5,1) circle (1pt);
\draw [fill=white] (1.5,1.5) circle (1pt);
\end{tikzpicture}
\centering{$d>4$}
\end{minipage}
\caption{Regions corresponding to (\ref{regionlwp}) for each case of dimension $d$.}\label{Region-d>=4Z}
\end{figure}
%
\end{center}

In the next table, we list the best results to date (as far as we know) for the system \eqref{Z}. 

\begin{table}[htp] 
\begin{tabular}{||c||rp{12.5cm}||}
\hline 
\cline{1-3} 
\cline{1-3} 
\multicolumn{1}{ ||c||  }
{\multirow{2}{*}{$d=1$}} 
& 
l.w.p.: 
& 
Theorem \ref{gtv1}.  
\vspace{0.15cm}
\\ 
\cline{2-3} 
\multicolumn{1}{||c||}{} 
& 
ill-p.:
& 
Biagioni and Linares proved in \cite{BL} non-existence of uniformly continuous \textit{solution mapping}, 
for $k<0$ and $l\le-3/2$; 
Holmer proved in \cite{holmer} norm inflation for $0<k<1$ and $l>2k-1/2$ and for $k\le0$ and $l>-1/2$; 
Also in \cite{holmer}, non-existence of uniformly continuous \textit{solution mapping} is proved for $k=0$ and $l<-3/2$; 
Theorem \ref{main02} (see Remark\ref{rmksenseill}) and Theorem \ref{main01}. 
\vspace{0.15cm}
\\
\hline 
\cline{1-3} 
\cline{1-3} 
\multicolumn{1}{ ||c||  }{\multirow{2}{*}{$d=2$}} 
& 
l.w.p.: 
& 
For $(k,l)\!=\!(0,-1/2)$, 
proved by 
Bejenaru, Herr, Holmer and Tataru in \cite{BHHT}; 
Theorem \ref{gtv1}. 
\vspace{0.15cm}
\\ 
\cline{2-3} 
\multicolumn{1}{||c||}{} 
& 
ill-p.:
& 
Theorem \ref{main02} (see Remark\ref{rmksenseill}) and Theorem \ref{main01}. 
\vspace{0.15cm}
\\
\hline 
\cline{1-3} 
\cline{1-3} 
\multicolumn{1}{ ||c||  }{\multirow{2}{*}{$d=3$}} 
& 
l.w.p.: 
& 
Theorem \ref{gtv1}.  
\vspace{0.15cm}
\\ 
\cline{2-3} 
\multicolumn{1}{||c||}{} 
& 
ill-p.:
& 
Theorems \ref{main02} and \ref{main01}. 
\vspace{0.15cm}
\\
\hline 
\cline{1-3} 
\cline{1-3} 
\multicolumn{1}{ ||c||  }{\multirow{2}{*}{$d=4$}} 
& 
l.w.p.: 
& 
For $l\ge0\,$, $k<4l+1$, $\max\{(l+1)/2\,,\,l-1\}\le k\le\min\{l+2,2l+11/8\}$ and $(k,l)\ne(2,3)$, 
proved by Bejenaru, Guo, Herr and Nakanishi in \cite{BGHN}; 
Theorem \ref{gtv1}. 
\vspace{0.15cm}
\\ 
\cline{2-3} 
\multicolumn{1}{||c||}{} 
& 
ill-p.:
& 
For $(k,l)=(2,3)$, 
non-existence of solution is proved in \cite{BGHN}; 
Theorems \ref{main02} and \ref{main01}. 
\vspace{0.15cm}
\\
\hline 
\cline{1-3} 
\cline{1-3} 
\multicolumn{1}{ ||c||  }{\multirow{2}{*}{$d>4$}} 
& 
l.w.p.: 
& 
Theorem \ref{gtv1}.  
\vspace{0.15cm}
\\ 
\cline{2-3} 
\multicolumn{1}{||c||}{} 
& 
ill-p.:
& 
Theorems \ref{main02} and \ref{main01}. 
\vspace{0.15cm}
\\ 
\cline{1-3} 
\cline{1-3} 
\hline 
\end{tabular}
\vspace{0.15cm}
\caption{Best l.w.p. (local well-posedness) and ill-p. (ill-posedness) results known to date for \eqref{Z}.}\label{bestresultstodate} 
\end{table}
	For $d\ge4$, Kato and Tsugawa in \cite{KS} proved the global well-posedness of the Zakharov system   
	for small data in 
	the mixed inhomogeneous and homogeneous space 
	$H^k(\R ^d)\!\times\!\dot H^l(\R ^d)\!\times\!\dot H^{l\!-\!1}\!(\R ^d)$ 
	at critical regularity $(k,l)\!=\!(\,d/2-3/2\,,\,d/2-2)$. 
	Global well-posedness for the Zakharov system is also studied in 
	\cite{pecher01}, \cite{pecher02}, \cite{cohotzi}, \cite{FPZ}, \cite{kishimoto01} and \cite{BGHN}. 

	Now we start to state our results. 
	First, we outline some definitions. 
	Assume that the system \eqref{Z} is 
	locally well-posed in the time interval $[0,T]$. 
	Then the \textit{solution mapping} associated to the system \eqref{Z} is the following map 
	\begin{eqnarray}
		\label{datasol}
		S\quad\!\! :\quad 
		B_r 
		\quad 
		&\longrightarrow& 
		\quad\mathcal C
		(\ [0,T]\ ;\ 
		H^{k,l}\ 
		) 
		\\ 
		\ \ \ \ \ \ (\varphi,\psi,\phi)\!\! 
		&\mapsto& 
		(u_{\!_{(\varphi,\psi,\phi)}}\,,\,n_{\!_{(\varphi,\psi,\phi)}}\,,\,\partial_tn_{\!_{(\varphi,\psi,\phi)}})\,,\nonumber 
	\end{eqnarray}
	where 
	$\mathcal C([0,T]\ ;\ H^{k,l})$ is a short notation for $C([0,T];H^k(\R^d))\!\times C([0,T];H^l(\R^d))\!\times C([0,T];H^{l-1}(\R^d))$, \linebreak 
	$B_r=\{(\varphi,\psi,\phi)\in H^{k,l}: \|(\varphi,\psi,\phi)\|_{H^{k,l}}<r\}$ 
	and 
	$u_{\!_{(\varphi,\psi,\phi)}}$ and $n_{\!_{(\varphi,\psi,\phi)}}$ are local solutions$^1$ 
\footnotetext[1]
{\ Precisely, $u_{\!_{(\varphi,\psi,\phi)}}$, $n_{\!_{(\varphi,\psi,\phi)}}$, $\partial_tn_{\!_{(\varphi,\psi,\phi)}}$satisfy the integral equations \eqref{uduhamel}, \eqref{nduhamel}, \eqref{ntduhamel} associated to the system \eqref{Z}, for all $t\in[0,T]$.}
	for system \eqref{Z} 
	with initial data $(u,v,\partial_tn)|_{t=0}=(\varphi,\psi,\phi)$. 
	
	Since Theorem \ref{gtv1} was obtained by means of contraction method, one can conclude the following: If $(k,l)$ satisfies conditions \eqref{regionlwp} then for every fixed $r>0$ there is a $T=T(r,k,l)>0$ such that the \textit{solution mapping} \eqref{datasol} is analitic (see Theorem. 3 in \cite{BTao}). 
So, if the system \eqref{Z} is locally well-posed in $H^{k,l}$ and the \textit{solution mapping} \eqref{datasol} fails to be $m$-times differentiable, then the usual contraction method can not be applied to prove the local well-posedness. In this case, we have a sense of ill-posedness and we say that the system \eqref{Z} is ill-posed by the method or simply the system \eqref{Z} is \textit{$C^m$ ill-posed}$^2$ 
\footnotetext[2]
{\ Actually, $C^m$ ill-posedness means that the \textit{solution mapping} is not $m$-times Fréchet differentiable.} 
in $H^{k,l}$. 

	Now fix $t\in[0,\,T]$. 
	Hereafter we call \textit{flow mapping} associated to the system \eqref{Z} the following map 
	\begin{eqnarray}
		\label{fluxo}
		S^{\,t}\!\quad\!\! :\quad 
		B_r
		\quad 
		&\longrightarrow& 
		\ \ \,H^k(\R^d)\times H^l(\R^d)\times H^{l-1}(\R^d)
		\\ 
		\ \ \ \ \ \ (\varphi,\psi,\phi)\!\! 
		&\mapsto& 
		\left(u_{\!_{(\varphi,\psi,\phi)}}(t)\,,\,n_{\!_{(\varphi,\psi,\phi)}}(t)\,,\,\partial_tn_{\!_{(\varphi,\psi,\phi)}}(t)\right).\nonumber 
	\end{eqnarray}
	We are now ready to enunciate our results. 
	Our first theorem shows that, in any dimension, 
	the regularity $(k,l)=(0,-1/2)$ is the endpoint for achieving well-posedness by contraction method 
	(see Figure \ref{Region-StnotC2}). 
	 \begin{theorem}\label{main02} 
	Let $d\in\N$. 
	Assume that the system \eqref{Z} is locally well-posed in the time interval $[0,T]$. 
	For any fixed $t\in (0,T]$, the flow mapping \eqref{fluxo} fails to
	 	be $C^2$ 
		at 
		the origin in $H^{k,l}$, provided 
		$l<-1/2\ $ or $\ l>2k-1/2\ $.	
	 \end{theorem}
	According to \cite{GTV97} (see p. 387), 
	the optimal relation between $k$ and $l$ is $l-k+1/2=0$. 
	Our next theorem shows that when $|l-k+1/2|>3/2$ (i.e., $l<k-2\ $ or $\ l>k+1$) the system \eqref{Z} is $C^2$ ill-posed (see Figure \ref{Region-SnotC2}). 
	\begin{theorem}\label{main01} 
	Let $d\in\N$. 
	Assume that the system \eqref{Z} is locally well-posed in the time interval $[0,T]$. 
	The solution mapping \eqref{datasol} fails to
		be $C^2$ 
		at 
		the origin in $H^{k,l}$, provided $l<k-2\ $ or $\ l>k+1$. 
	\end{theorem}
\begin{remark}\label{rmksenseill}
The sense of ill-posedness stated in Theorem \ref{main02} is slightly stronger than the sense 
stated in Theorem \ref{main01}. Indeed, if the \textit{flow mapping} \eqref{fluxo} is not $C^2$, 
neither is, a fortiori, 
the \textit{solution mapping} \eqref{datasol}. 
Thus, Theorem \ref{main02} slightly improves the ill-posedness results in \cite{holmer} and \cite{BHHT}, 
for $d=1$ and $d=2$, respectively, 
both establishing that the \textit{solution mapping} \eqref{datasol} is not $C^2$ for $l<-1/2\ $ or $\ l>2k-1/2\ $. 
\end{remark}
\begin{remark}
Theorem \ref{main01} stablishes $C^2$ ill-posedness for new indices $(k,l)$ (see Figure \ref{Region-SnotC2}). 
For such indices, the difference of regularity between the initial data is large 
(i.e., $l\gg k$ or $k\gg l$). 
Such result seems natural, due to coupling of the system via nonlinearities. 
Indeed, 
for instance, 
high regularity for $u(t)$ is not expect when $n(t)$ has low regularity, 
in view of \eqref{uduhamel}. 
By the way, the $C^2$ ill-posedness for $l<k-2$ is obtained by dealing with \eqref{uduhamel}. 
\end{remark}
\begin{remark}
In the periodic setting, 
Kishimoto proved in \cite{kishimoto02} 
the 
$C^2$ ill-posedness$^3$ 
of the Zakharov system
in $H^k\!(\mathbb T^d)\!\!\times\!H^l\!(\mathbb T^d)\!\!\times\!H^{l\!-\!1}\!(\mathbb T^d)$ for $d\ge2$, 
provided $\ l<\max\{0\,,\,k-2\}\ $ or $\ l>\min\{2k-1\,,\,k+1\}$. 
These indices $(k,l)$ are exactly the same of Theorems \ref{main02} and \ref{main01}, 
excepting for admitting $-1/2\le l < 0$. 
We point out that 
in \cite{BHHT} was proved, by means of contraction method, that the system \eqref{Z} is locally well-posed for $d=2$, $k=0$ and $l=-1/2$. 
\end{remark}
\footnotetext[3]
{\ $C^2$ ill-posedness in the slightly weaker sense (see Remark \ref{rmksenseill}). 
However, 
for $d=2$ and particular $(k,l)$ is proved in \cite{kishimoto02} 
ill-posedness in much stronger senses, namely 
norm inflation and non-existence of continuous \textit{solution mapping}.} 

This paper is organized as follows. 
In Section \ref{secnotations}, we introduce some notations to be used throughout the whole text. 
In Section \ref{secprelmanal}, is presented a preliminary analysis which provides a methodical approach to our proofs, exposing the main ideas. 
In Section \ref{secproofmain02*}, we prove Theorem \ref{main02} and in Section \ref{secproofmain01}, we prove Theorem \ref{main01}.

	\section{Notations}\label{secnotations}
 	%
%
%
\begin{itemize}
\item 
$(*.*)_R\ $(or $(*.*)_L$) denotes the right(or left)-hand side of an equality or inequality numbered by $(*.*)$. \vspace{0.2cm} 

\item 
$\|(\varphi,\psi,\phi)\|_{H^{k,l}}^2=\|\varphi\|_{H^{k}}^2+\|\psi\|_{H^{l}}^2+\|\phi\|_{H^{l-1}}^2$, 
where $H^{k,l}\, =\, H^k\!(\R ^d;\mathbb C)\times H^l\!(\R ^d;\R)\times H^{l\!-\!1}\!(\R ^d;\R)$. 
\vspace{0.2cm} 

\item 
$\cjap{\xi}=\sqrt{1+|\xi|^2}\ $,\, $\xi\in\R^d$. 
\vspace{0.2cm} 

\item 
$\chi_{_{\Omega}}$ denotes the characteristic function of $\Omega\subset \R^d$. 
\vspace{0.2cm} 

\item 
$|\Omega|$ denotes de Lebesgue measure of the set $\Omega$, i.e., $|\Omega|=\int\chi_{_{\Omega}}(\xi)d\xi$. 
\vspace{0.2cm} 

\item 
$\mathcal S(\R^d)$ denotes the Schwartz space and $\mathcal S'(\R^d)$ denotes the space of tempered distributions.
\vspace{0.2cm} 

\item 
$\widehat{f}$ and $\check{f}$ denote, respectively, the Fourier transform and the inverse Fourier transform of $f\in\mathcal S'(\R^d)$. 
\end{itemize}

%

	\section{Preliminary Analysis}\label{secprelmanal}

	The integral equations associated to the system \eqref{Z} 
	with initial data $(u,v,\partial_tn)|_{t=0}=(\varphi,\psi,\phi)$ 
	are 
	\begin{eqnarray}
		u(t)&=&e^{it\Delta }\varphi\ -\ i\!\!\int_0^te^{i(t-s)\Delta}u(s)n(s)ds,\label{uduhamel}\\ 
		n(t)&=&W(t)(\psi,\phi)\ +\ \int_0^t
		W_1(t-s) 
		\Delta|u|^2(s)ds,\label{nduhamel}\\ 
		\partial_tn(t)&=&
		W(t)(\phi,\Delta\psi) 
		\ +\ 
		\int_0^t
		W_0(t-s) 		
		\Delta|u|^2(s)ds,\label{ntduhamel}
	\end{eqnarray}
	where $\{e^{it\Delta}\}_{t\in\R}$ is the unitary group 
	in $H^s(\R^d)$ 
	associated to the linear Schrödinger equation, given by \linebreak 
	$e^{it\Delta}\varphi:=\{e^{-it|\cdot|^2}\widehat{\varphi}(\cdot)\}\!\,\check{\ }$ 
	and 
	$\{W(t)\}_{t\in\R}$ is the linear wave propagator 
	$W(t)(\psi,\phi):=W_0(t)\psi+W_1(t)\phi$, 
	where $W_0$ and $W_1$ are given by $W_0(t)\psi=\cos\left(t\sqrt{-\Delta}\right)\psi:=\{\cos(t|\cdot|)\widehat{\psi}(\cdot)\}\!\,\check{\ }\ $ and $\ W_1(t)\phi=\tfrac{\sin\left(t\sqrt{-\Delta}\right)}{\sqrt{-\Delta}}\phi:=\left\{\tfrac{\sin(t|\cdot|)}{|\cdot|}\,\widehat{\phi}(\cdot)\right\}\!\check{\ }$. 

	Assume that the system \eqref{Z} is locally well-posed in $H^{k,l}$, in the time interval $[0,T]$. 
	Suppose also that there exists $t\in[0,T]$ such that the \textit{flow mapping} \eqref{fluxo} is two times Fréchet differentiable at the origin in $H^{k,l}$. 
	Then, the second Fréchet derivative of $S^{\,t}$ at origin belongs to $\mathcal B$, the normed space of bounded bilinear applications from $H^{k,l}\times H^{k,l}$ to $H^{k,l}$. 
	In particular, we have the following estimate for the second Gâteaux derivative of $S^{\,t}\!$ at origin 
\begin{equation}\label{ill18**}
	\left\| \frac{\partial S^{\,t}_{(0,0,0)}}{\partial\Phi_0\partial\Phi_1}\right\|_{_{H^{k,l}}}\!\! 
	=
	\left\| D^2S^{\,t}_{(0,0,0)}(\Phi_0,\Phi_1)\right\|_{_{H^{k,l}}}\!\! 
	\leq\ 
	\left\|  D^2S^{\,t}_{(0,0,0)}  \right\|_{\mathcal B} \|\Phi_0\|_{_{H^{k,l}}}\|\Phi_1\|_{_{H^{k,l}}}\,,
	\quad \forall\Phi_0,\Phi_1\in H^{k,l}.
\end{equation}
	Similarly, assuming \textit{solution mapping} \eqref{datasol} two times Fréchet differentiable at the origin, 
	we have $D^2S_{(0,0,0)}$ belonging to $\mathcal B_{\mathcal C}$, 
	the normed space of bounded bilinear applications from 
	$H^{k,l}\times H^{k,l}$ to 
	$\mathcal C([0,T];H^{k,l})$. 
	Then 
\begin{equation}\label{ill01**}
	\sup_{t\in [0,T]}\left\| \frac{\partial S^{\,t}_{(0,0,0)}}{\partial\Phi_0\partial\Phi_1}\right\|_{_{H^{k,l}}}\!\! 
	\leq\ 
	\left\|  D^2S_{(0,0,0)}  \right\|_{\mathcal B_{\mathcal C}} \|\Phi_0\|_{_{H^{k,l}}}\|\Phi_1\|_{_{H^{k,l}}}\ ,
	\quad\quad \forall\Phi_0,\Phi_1\in H^{k,l}.
\end{equation}

	Thus, we can prove Theorem \ref{main02} by showing that estimate \eqref{ill18**} is false for ${(k,l)}$ in the region of Figure \ref{Region-StnotC2}. 
	In the case of Theorem \ref{main01}, the indices $(k,l)$ in the region of Figure \ref{Region-SnotC2} impose additional technical difficulties to get good lower bounds for $\eqref{ill18**}_L$. 
	To overcome such difficulties, we made use of a sequence $t_N\to0$, in consequence, we merely prove that estimate \eqref{ill01**} is false, obtaining an ill-posedness result in a slightly weaker sense. 
	
	Since $S^{\,t}_{(0,0,0)}=(0,0,0)$, 
	for each direction $\Phi=(\varphi,\psi,\phi)\in\mathcal S(\R^d)\times\mathcal S(\R^d)\times\mathcal S(\R^d)$, 
	the first Gâteaux derivatives of $\eqref{uduhamel}_R$, $\eqref{nduhamel}_R$ and $\eqref{ntduhamel}_R$ at the origin are 
	$e^{it\Delta}\varphi\ $, 
	$W(t)(\psi,\phi)\,$ 
	and 
	$\,W(t)(\phi,\Delta\psi)$, 
	respectively. 
	Further, 
	from \eqref{ill18**}, 
	we deduce the following estimates 
	for the second Gâteaux derivatives of $u(t)$, $n(t)$ and $\partial_tn(t)$ 
	in the directions 
	$(\Phi_0,\Phi_1)\!=\!(\,(\varphi_0,\psi_0,\phi_0)\,,\,(\varphi_1,\psi_1,\phi_1)\,)\in(\mathcal S(\R^d)\times \mathcal S(\R^d)\times \mathcal S(\R^d))^2$ 
	\begin{eqnarray}
		\left\|\frac{\partial^2 u_{\!_{(0,0,0)}}}{\partial\Phi_0\partial\Phi_1}(t)\! \right\|_{_{\!H^k}}
		\!\!\!\!\!\!\!\!&=&\!\!\!\! 
		\left\|  \int_0^te^{i(t-s)\Delta}\{e^{is\Delta}\varphi_0W(s)(\psi_1,\phi_1)+e^{is\Delta}\varphi_1W(s)(\psi_0,\phi_0)\}ds   \right\|_{_{\!H^k}} 
		\!\lesssim 
		\|\Phi_0\|_{_{\!H^{k,l}}}\|\Phi_1\|_{_{\!H^{k,l}}},\ \ \label{ill03**}
		\\ 
		\left\|\frac{\partial^2 n_{\!_{(0,0,0)}}}{\partial\Phi_0\partial\Phi_1}(t)\! \right\|_{_{\!H^l}}
		\!\!\!\!\!\!\!\!&=&\!\!\!\! 
		\left\|\int_0^t
		W_1(t-s) 
		\Delta\{e^{is\Delta}\varphi_0\,\overline{e^{is\Delta}\varphi_1}+\overline{e^{is\Delta}\varphi_0}\,e^{is\Delta} \varphi_1\}ds   \right\|_{_{\!H^l}}
		\!\lesssim 
		\|\Phi_0\|_{_{\!H^{k,l}}}\|\Phi_1\|_{_{\!H^{k,l}}}, \label{ill04**}
		\\ 
		\left\|\frac{\partial^2\partial_tn_{\!_{(0,0,0)}}}{\partial\Phi_0\partial\Phi_1}(t)\! \right\|_{_{\!H^{l-1}}}
		\!\!\!\!\!\!\!\!\!\!\!\!&=&\!\!\!\! 
		\left\|\int_0^t
		W_0(t-s) 
		\Delta\{e^{is\Delta}\varphi_0\,\overline{e^{is\Delta}\varphi_1}+\overline{e^{is\Delta}\varphi_0}\,e^{is\Delta} \varphi_1\}ds   \right\|_{_{\!H^{l-1}}}
		\!\lesssim 
		\|\Phi_0\|_{_{\!H^{k,l}}}\|\Phi_1\|_{_{\!H^{k,l}}}. \label{ill500}
	\end{eqnarray} 

	Hence, the proof of Theorem \ref{main02} boils down to 
	getting sequences of directions $\Phi$ showing that one of these last three estimates fails for the fixed $t\in[0,T]$. 
	For Theorem \ref{main01}, such sequences just need to show that one of \eqref{ill03**}-\eqref{ill500} can not hold uniformly for $t\in[0,T]$. 
	
	We deal with \eqref{ill03**} by choosing directions $\Phi_0=\Phi_1=(\varphi,\psi,0)$ with $\varphi,\psi\in S(\R^d)$. 
	Since in $\mathcal S(\R^d)$ the Fourier transform convert products in convolutions, 
	from \eqref{ill03**} we conclude the following estimate 
\begin{equation}\label{illconvschwarz1} 
	\left\|\!\cjap{\xi}^k\!\!\!\int_0^t\!\!e^{-i(t-s)|\xi|^2}\!\!\!\int_{\!\R^d}e^{-is|\xi_1|^2}\widehat{\varphi}(\xi_1)\cos(s|\xi-\xi_1|)\widehat{\psi}(\xi-\xi_1)d\xi_1ds     \right\|_{_{\!L^2_\xi}} 
	\!\lesssim\ 
	\left\|\varphi\right\|_{_{\!H^k}}^2\!\!+\!\left\|\psi\right\|_{_{\!H^l}}^2, 
	\,\ \forall\varphi,\psi\in\mathcal S(\R^d). 
\end{equation}
Hereafter we will denote, as usual, $\xi_2\ :=\ \xi-\xi_1$, then 
\begin{equation}\label{xi1+xi2=xi} 
\xi_1+\ \xi_2\ =\ \xi\,.
\end{equation} 

	For bounded subsets 
	$A,B\subset\R^d$, by taking $\varphi,\psi\in\mathcal S(\R^d)$ such that$^4$ 
	$\langle\cdot\rangle^k\,\widehat{\varphi}\ \sim\ \chi_{_{A}}$
	and $\langle\cdot\rangle^l\,\widehat{\psi}\ \sim\ \chi_{_{B}}$, 
	we conclude from \eqref{illconvschwarz1} that 
%
%
\footnotetext[4]
{\ Precisely, 
$\chi_{_{A}}\leq\langle\cdot\rangle^k\,\widehat{\varphi}$\, with\, $\|\varphi\|_{H^k}\leq2\|\chi_{_{A}}\|_{L^2}$ 
\ and\ 
$\chi_{_{B}}\leq\langle\cdot\rangle^l\,\widehat{\psi}$\, with\, $\|\psi\|_{H^l}\leq2\|\chi_{_{B}}\|_{L^2}$.} 
\begin{equation}\label{ill05**}
	\left\| \int_0^t\int_{\!\R^d}\tfrac{\cjap{\xi}^k}{\cjap{\xi_1}^k\cjap{\xi_2}^l}\cos(s|\xi|^2\!-\!s|\xi_1|^2)\cos(s|\xi_2|)\chi_{_{A}}(\xi_1)\chi_{_{B}}(\xi_2)d\xi_1ds\right\|_{_{\!L^2_\xi}} 
	\!\lesssim\ 
	|A|+|B|\,.
\end{equation} 

We can rewrite $\eqref{ill05**}_L$ as 
\begin{equation}\label{cosprodsumcos*} 
	\left\| 
	\int_0^t\int_{\!\R^d}\tfrac{\cjap{\xi}^k}{\cjap{\xi_1}^k\cjap{\xi_2}^l}\, 
	\tfrac12 
	\left[
	\cos(\sigma_+s)+\cos(\sigma_-s)
	\right]\,
	\chi_{_{A}}(\xi_1)\chi_{_{B}}(\xi_2)d\xi_1ds\right\|_{_{\!L^2_\xi}}, 
\end{equation} 
where $\sigma_+$ and $\sigma_+$ are what we call \textit{the algebraic relations} associated to \eqref{ill03**}, given by 
\begin{equation}\label{sigma+-1*} 
\sigma_\pm\ :=\ |\xi|^2-|\xi_1|^2\pm|\xi_2|. 
\end{equation} 

	Finally, we have to choose sequences of sets $\{A_N\}_{N\in\N}$ and $\{B_N\}_{N\in\N}$ 
	such that, for $\xi_1\in A_N$ and $\xi_2\in B_N$, yields 
	increasing $\tfrac{\cjap{\xi}^k}{\cjap{\xi_1}^k\cjap{\xi_2}^l}$, 
	small $\sigma_+$ and large $\sigma_-$, when $N\to+\infty$. 
	It allows us to get good lower bounds for \eqref{cosprodsumcos*}, since 
\begin{equation}\label{cossincontrol*} 
\cos(\theta)>1/2\,,
\quad  
\forall\theta\in(-1,1) 
\quad\quad\quad 
\textrm{and}\ 
\quad\quad\quad 
\int_0^t\cos(ks)ds=\frac{\sin(kt)}k\,,
\quad 
\forall k\ne0. 
\end{equation}
	Moreover, we will need a lower bound for $\|\chi_{_{_{\!A_N}}}\!\!\ast\,\chi_{_{_{\!B_N}}}\|_{L^2}$. 
	For this purpose, the next elementary result is very useful.
\begin{lemma}\label{lemaleandro*}
	(\cite{leandro}) Let $A,B,R\subset \R^d$. If $R-B=\{x-y\ :\  x\in R\ \ \textrm{and}\ \ y\in B\}\subset A$ then 
		$$|R|^{\frac12}|B|\leq \|\chi_{_{\!A}}\ast\chi_{_{\!  B}}\|_{_{\!L^2(\R^d)}}.$$	
\end{lemma}
\begin{remark}\label{rmkfimprelimanal}
	For the case $l<-1/2$ in Theorem \ref{main02}, 
	by a good choice of $A_N$ and $B_N$, 
	it is possible to obtain a \lq\lq high + high = high\rq\rq\ interaction in \eqref{xi1+xi2=xi} 
	providing  
	\lq\lq high\rq\rq\,  $\tfrac{\cjap{\xi}^k}{\cjap{\xi_1}^k\cjap{\xi_2}^l}$\,, 
	\lq\lq low\rq\rq\, $\sigma_+$ and \lq\lq high\rq\rq\, $\sigma_-$, 
	which yield good lower bounds for \eqref{cosprodsumcos*}. 
	But for the case $k-l>2$ in Theorem \ref{main01}, 
	to obtain 
	\lq\lq high\rq\rq\, $\tfrac{\cjap{\xi}^k}{\cjap{\xi_1}^k\cjap{\xi_2}^l}$, 
	the interaction 
	must be of type \lq\lq low + high = high\rq\rq, 
	implying 
	\lq\lq high\rq\rq\, $\sigma_+$ and \lq\lq high\rq\rq\, $\sigma_-$,
	which do not provide lower bound for \eqref{cosprodsumcos*}. 
	Then we choose a sequence $t_N\to0$, 
	allowing us to obtain lower bounds directly from $\eqref{ill05**}_L$. 
\end{remark} 
%

	\section{Proof of Theorem \ref{main02}}\label{secproofmain02*}

	Assume that, for a fixed $t\in(0,T]$, the \textit{flow mapping} \eqref{fluxo} is $C^2$ at the origin. 
	Then, from \eqref{ill05**}, \eqref{cosprodsumcos*} and \eqref{sigma+-1*}, we get the following estimate for bounded subsets $A,B\subset\R^d$
\begin{equation}\label{ill1*}
	\|I_{_{A\,,\,B}}^+(\xi)\|_{_{\!L^2_\xi}}\ -\ \|I_{_{A\,,\,B}}^-(\xi)\|_{_{\!L^2_\xi}} 
	\ \ \lesssim\ \ 
	|A|+|B|\,, 
\end{equation}
where 
%
\begin{equation}\label{ill19*}
I_{_{A\,,\,B}}^\pm(\xi) 
\ :=\ 
\int_0^t\int_{\!\R^d}\tfrac{\cjap{\xi}^k}{\cjap{\xi_1}^k\cjap{\xi_2}^l}\cos(\sigma_\pm s)\chi_{_{A}}(\xi_1)\chi_{_{B}}(\xi_2)d\xi_1ds. 
\end{equation} 

Note that, for $\xi_1=(\xi_1^1,\cdots,\xi_1^d)\in \R^d$ and $\ \xi_2=(\xi_2^1,\cdots,\xi_2^d)\in \R^d$, we can rewrite \eqref{sigma+-1*} as  
\begin{equation}\label{sigmaaberto*}
\sigma_\pm
\ \ =\ \ 
\sum_{j=1}^d\!\!\left(|\xi_1^j+\xi_2^j|^2-|\xi_1^j|^2\right)\ \pm\ |\xi_2| 
\ \ =\ \ 
\xi_2^1(2\xi_1^1+\xi_2^1\pm1)
\ \pm\ 
(|\xi_2|-\xi_2^1)
\ +\ 
\sum_{j=2}^d\xi_2^j(2\xi_1^j+\xi_2^j). 
\end{equation} 

In order to obtain a lower bound for $\|I_{_{A\,,\,B}}^+\|_{_{\!L^2}}$ 
and an upper bound $\|I_{_{A\,,\,B}}^-\|_{_{\!L^2}}$, 
we choose the sets $A,B\subset\R^d$ taking \eqref{sigmaaberto*} into account. 
So, for $N\in\N$ and $0<\delta<\min\{\frac1{7t},1\}$, we define$^5$ 
\footnotetext[5]
{\ Evidently, if $d=1$ then $A$ and $B$ are just intervals, the last sum in \eqref{sigmaaberto*} does not exist and $\eqref{res1resto*}_R$ should be ignored.} 
\begin{equation*}
	A=A_N\ :=\ \left[-N\!\ ,\!\ -N+\tfrac{\delta}N\right]\times\left[0\!\ ,\!\ \tfrac{\delta}{d-1}\right]^{d-1} 
\hspace{0.25cm}\textrm{and}\hspace{0.5cm}
	B=B_N\ :=\ \left[2N-1\!\ ,\!\ 2N-1+\tfrac{\delta}{2N}\right]\times\left[0\!\ ,\!\ \tfrac{\delta}{2(d-1)}\right]^{d-1}.
\end{equation*}
Then, 
for $(\xi_1,\xi_2)\in A_N\times B_N$, 
we have 
\begin{equation}\label{ordfreq1*}
\cjap{\xi_1}\ \sim\ \cjap{\xi_2}\ \sim\ \cjap{\xi_1+\xi_2}\ \sim\ N
\end{equation}
and since $\delta<1$ we also have $\xi_2^1\in[N\!\ ,\!\ 2N]$ and $(2\xi_1^1+\xi_2^1)\in[-1\!\ ,\!\ -1+\frac{5\delta}{2N}]$. 
Thus, 
\begin{equation}\label{resson+-1*}
\xi_2^1(2\xi_1^1+\xi_2^1+1)
\in\ 
[0\!\ ,\!\ 5\delta]\ , 
\hspace{2cm} 
\xi_2^1(2\xi_1^1+\xi_2^1-1)
\in\ 
[-4N\!\ ,\!\ -N]\ ,
\end{equation}
\begin{equation}\label{res1resto*}
(|\xi_2|-\xi_2^1)\in\left[0\!\ ,\!\ \tfrac{\delta}2\right]
\hspace{0.75cm}\textrm{and}\hspace{0.75cm}
\sum_{j=2}^d\xi_2^j(2\xi_1^j+\xi_2^j)\in\left[0\!\ ,\!\ \tfrac{5\delta^2}{4(d-1)}\right].
\end{equation}
Therefore, combining \eqref{sigmaaberto*}, $\eqref{resson+-1*}_L$ and \eqref{res1resto*} we obtain  
\begin{equation}\label{sigm+s,1*}
\sigma_+\in[0\!\ ,\!\ 7\delta)
\end{equation}
and combining \eqref{sigmaaberto*}, $\eqref{resson+-1*}_R$ and \eqref{res1resto*} we obtain 
\begin{equation}\label{sigm-s,1*}
\sigma_-\in\left(-5N\!\ ,\!\ -\tfrac 12N\right). 
\end{equation}

Since $\delta<\frac1{7t}$, from \eqref{ill19*}, $\eqref{cossincontrol*}_L$, \eqref{sigm+s,1*} and \eqref{ordfreq1*} we get that 
\begin{equation}\label{I+lowboun*}
I_{_{A\,,\,B}}^+(\xi)\ \ge\ \tfrac12\int_0^t\int_{\!\R^d}\tfrac{\cjap{\xi}^k}{\cjap{\xi_1}^k\cjap{\xi_2}^l}\chi_{_{A}}(\xi_1)\chi_{_{B}}(\xi_2)d\xi_1ds
\ \gtrsim\ tN^{-l}\chi_{_{\!A}}\ast\chi_{_{\!B}}(\xi).
\end{equation}

Now, Lemma \ref{lemaleandro*} allows us to get a lower bound for $I_{_{A\,,\,B}}^+(\xi)$. For this purpose, consider the set 
\begin{equation*}
	R=R_N\ :=\ \left[N-1+\tfrac{\delta}{2N}\!\ ,\!\ N-1+\tfrac{\delta}{N}\right]\times\left[\tfrac{\delta}{2(d-1)}\!\ ,\!\ \tfrac{\delta}{d-1}\right]^{d-1}.
\end{equation*}
Then we have $R-B\subset A$ and 
\begin{equation}\label{mesureRAB1*} 
|R|\sim|A|\sim |B|\sim N^{-1}. 
\end{equation} 

Using \eqref{I+lowboun*}, Lemma \ref{lemaleandro*} and \eqref{mesureRAB1*} we obtain that
\begin{equation}\label{ill23*}
	\|I_{_{A\,,\,B}}^+\|_{_{L^2}}\ \gtrsim\ tN^{-l}|R|^{\frac12}|B|\sim\ t N^{-l-\frac{3}{2}}.
\end{equation}

On the other hand, using 
\eqref{ill19*}, 
the Fubini's theorem, $\eqref{cossincontrol*}_R$, \eqref{ordfreq1*}, \eqref{sigm-s,1*}, Young's convolution inequality and \eqref{mesureRAB1*}, we get that 
\begin{equation}\label{ill26*}
	\|I_{_{A\,,\,B}}^-\|_{_{L^2}}
	\ =\ 
	\left\|
	\int_{\R^d}\tfrac{\cjap{\xi}^k}{\cjap{\xi_1}^k\cjap{\xi_2}^l}\tfrac{\sin(\sigma_-t)}{\sigma_-}\chi_{_{A}}(\xi_1)\chi_{_{B}}(\xi_2)d\xi_1\right\|_{_{L^2_\xi}}
	\ \lesssim\ 
	\left\|\frac1{N^l}\frac1{N^{\ }}\ \chi_{_{\!A}}\ast\chi_{_{\!B}}\right\|_{_{L^2}}
         \ \le\ 
         \frac{|A||B|^{\frac12}}{N^{l+1}}
         \ \sim\  
         N^{-l-\frac52}.
\end{equation}

Finally, combining \eqref{ill1*}, \eqref{ill23*}, \eqref{ill26*} and \eqref{mesureRAB1*} we conclude that
\begin{equation*}
	tN^{-l-\frac32}-N^{-l-\frac52}\ \lesssim\ 
	N^{-1}, 
	\quad\quad\forall N\in\N.
\end{equation*}
Hence $l\ \ge\ -1/2$ when the \textit{flow mapping} \eqref{fluxo} is $C^2$ at the origin.
\\ 

Now we will show that $l\le2k-1/2$ dealing with \eqref{ill04**}. 
	Similarly to the manner that we obtained \eqref{illconvschwarz1}, 
	using now $\Phi_0=(\varphi,0,0)$ and $\Phi_1=(\upsilon,0,0)$ in \eqref{ill04**} 
	with $\varphi,\upsilon\in\mathcal S(\R^d)$, we obtain  
\begin{equation*}\label{illconvschwarz2} 
	\left\|  \! \cjap{\xi}\!^l\!\!\!\int_0^t\!
	\tfrac{e^{i(t-s)|\xi|}-e^{\!-i(t-s)|\xi|}}{2i|\xi|}|\xi|^2\!\!\!
	\int_{_{\!\!\R^d}}\!\!\!\!e^{-is|\xi_1|^2}\widehat{\varphi}(\xi_1)e^{is|\xi_2|^2}\overline{\widehat{\upsilon}(-\xi_2)}+e^{is|\xi_1|^2}\overline{\widehat{\varphi}(-\xi_1)}e^{-is|\xi_2|^2}\widehat{\upsilon}(\xi_2)d\xi_1ds \right\|_{_{\!L^2_\xi}}
	\!\!\lesssim 
	\left\|\varphi\right\|_{_{\!H^k}}\left\|\upsilon\right\|_{_{\!H^l}}. 
\end{equation*}
Similarly to \eqref{illconvschwarz1} and \eqref{ill05**}, from the last estimate follows that, 
for bounded subsets $A,B\subset\R^d$, we have 
\begin{equation*}\label{ill***} 
	\Biggl\|\!\int_0^t\!\!\!\int\!\!\!\tfrac{\cjap{\xi}^l|\xi|}{\cjap{\xi_1}^k\cjap{\xi_2}^k}\!\!
	\left(
	\!e^{i(t-s)|\xi|}\!-\!e^{\!-i(t-s)|\xi|}\!
	\right)\!\!\!
	\left(
	\!e^{\!-is(|\xi_1|^2\!-\!|\xi_2|^2)}\chi_{_{_{\!A}}}\!\!(\xi_1)\chi_{_{_{\!-B}}}\!\!(\xi_2)+e^{is(|\xi_1|^2\!-\!|\xi_2|^2)}\chi_{_{_{\!-A}}}\!\!(\xi_1)\chi_{_{_{\!B}}}\!\!(\xi_2)\!
	\right)\!\!
	d\xi_1ds\Biggr\|_{_{\!L^2_\xi}}\!\! 
	\lesssim  
	|A|^{\frac12}|B|^{\frac12}\,. 
	\end{equation*} 
So, under the additional assumption that the sets $(A+(-B))$ and $((-A)+B)$ are disjoint, we have$^6$ 
\footnotetext[6]
{\ Since\ 
$\chi_{_{_{\!X}}}\!\!\ast\,\chi_{_{_{\!Y}}}
=
\chi_{_{_{\!X+Y}}}\,\chi_{_{_{\!X}}}\!\!\ast\,\chi_{_{_{\!Y}}}$\ 
and\ 
$\|f\,\chi_{_{_{\!Z}}}+g\,\chi_{_{_{\!W}}}\|_{L^2}^2=\|f\,\chi_{_{_{\!Z}}}\|_{L^2}^2+\|g\,\chi_{_{_{\!W}}}\|_{L^2}^2$\ 
when\ 
$Z\cap W=\emptyset$
.} 
\begin{equation}\label{ill2*}
	\|J^+_{_{A,B}}(\xi)\|_{_{\!L^2_\xi}}-\|J^-_{_{A,B}}(\xi)\|_{_{\!L^2_\xi}}\ 
\le\ 
\Biggl\|\int_0^t\!\!\!\int_{_{\!\!\R^d}}\!\!\tfrac{\cjap{\xi}^l|\xi|}{\cjap{\xi_1}^k\cjap{\xi_2}^k}
	\!\left(
	\!e^{it|\xi|-is\zeta_+}\!-\!e^{-it|\xi|-is\zeta_-}
	\right)\!
	\chi_{_{_{\!A}}}\!\!(\xi_1)\chi_{_{_{\!-B}}}\!\!(\xi_2)
	d\xi_1ds\Biggr\|_{_{\!L^2_\xi}} 
\lesssim  
	|A|^{\frac12}|B|^{\frac12}\,, 
\end{equation}
where $\zeta_+$ and $\zeta_-$ are \textit{the algebraic relations} associated to \eqref{ill04**} given by 
\begin{equation}\label{zeta+-*}
\zeta_\pm\ :=\ |\xi_1|^2-|\xi_2|^2\pm|\xi|
\ =\ \xi^1(\xi_1^1-\xi_2^1\pm1)
\ \pm\ 
(|\xi|-\xi^1)
\ +\ 
\sum_{j=2}^d\xi^j(\xi_1^j-\xi_2^j) 
\end{equation}
and 
%
\begin{equation*}\label{J+*}
J^\pm_{_{A,B}}(\xi)
\ :=\ |\xi|\!\int_0^t\int_{\!\R^d}\tfrac{\cjap{\xi}^l}{\cjap{\xi_1}^k\cjap{\xi_2}^k}e^{-is\zeta_\pm}
\chi_{_{_{\!A}}}\!\!(\xi_1)\chi_{_{_{\!-B}}}\!\!(\xi_2)d\xi_1ds. 
\end{equation*}

Now, in view of \eqref{zeta+-*}, we choose the sets $A$ and $B$. So, for $N\in\N$ and $0<\delta<\min\{\frac1{7t},1\}$, we define 
\begin{equation*}
	A=A_N\ :=\ \left[N\!\ ,\!\ N+\tfrac{\delta}N\right]\times\left[0\!\ ,\!\ \tfrac{\delta}{d-1}\right]^{d-1} 
\hspace{0.25cm}\textrm{and}\hspace{0.5cm}
	B=B_N\ :=\ \left[-N-1\!\ ,\!\ -N-1+\tfrac{\delta}{2N}\right]\times\left[-\tfrac{\delta}{2(d-1)}\!\ ,\!\ 0\right]^{d-1}. 
\end{equation*}
Then $(A+(-B))\cap((-A)+B)=\emptyset\ $ and $\cjap{\xi_1}\sim\cjap{\xi_2}\sim\cjap{\xi_1+\xi_2}\sim N$, 
for $(\xi_1,\xi_2)\in A_N\times B_N$.
Moreover, following the procedure used in \eqref{sigmaaberto*}-\eqref{sigm-s,1*}, one can verify that 
$\zeta_+\in(-\delta\!\ ,\!\ 7\delta)$ and $\zeta_-\in(-7N\!\ ,\!\ -N)$. 
Therefore, we have 
\begin{equation}\label{J+lowboun*}
|J^+_{_{A,B}}(\xi)| 
\ \gtrsim\ tN^{l-2k+1}\chi_{_{\!A}}\ast\chi_{_{\!B}}(\xi).
\end{equation}

Consider the set 
\begin{equation*}
	R=R_N\ :=\ \left[2N+1\!\ ,\!\ 2N+1+\tfrac{\delta}{2N}\right]\times\left[\tfrac{\delta}{2(d-1)}\!\ ,\!\ \tfrac{\delta}{(d-1)}\right]^{d-1}
\end{equation*}
and note that $R-(-B)\subset A$ and $|R|\sim|A|\sim |B|\sim N^{-1}$. 
Then, using \eqref{J+lowboun*} and Lemma \ref{lemaleandro*}, we obtain that
\begin{equation}\label{ill40*}
	\|J^+_{_{A,B}}\|_{_{L^2}}\ \gtrsim\ tN^{l-2k+1}|R|^{\frac12}|B|\sim\ t N^{l-2k-\frac12}.
\end{equation}

On the other hand, similarly to \eqref{ill26*}, we get that 
\begin{equation}\label{ill41*}
	\|J^-_{_{A,B}}\|_{_{L^2}}
	\ =\ 
	\left\||\xi|\int_{\R^d}\tfrac{\cjap{\xi}^l}{\cjap{\xi_1}^k\cjap{\xi_2}^k}\tfrac{(e^{-it\zeta_-}-1)}{-i\zeta_-}\chi_{_{_{\!A}}}\!\!(\xi_1)\chi_{_{_{\!-B}}}\!\!(\xi_2)d\xi_1\right\|_{_{L^2_\xi}}
	\ \lesssim\ 
         N^{l-2k-\frac32}.
\end{equation}

Finally, combining \eqref{ill2*}, \eqref{ill40*} and \eqref{ill41*} we conclude that
\begin{equation*}
	tN^{l-2k-\frac12}-N^{l-2k-\frac32}\ \lesssim\ |A|^{\frac12}|B|^{\frac12}\ \sim\ N^{-1}, 
	\quad\quad\forall N\in\N.
\end{equation*}
Hence $l\le2k-1/2$ when the \textit{flow mapping} \eqref{fluxo} is $C^2$ at the origin. \hfill $\blacksquare$

%
%
	\section{Proof of Theorem \ref{main01}}\label{secproofmain01}
	%

	Assume that the \textit{solution mapping} \eqref{datasol} is $C^2$ at the origin. 
	Employing the same procedure that yields \eqref{ill05**} from \eqref{ill18**}, 
	one can conclude, from \eqref{ill01**}, the following estimate 
	for bounded subsets $A,B\subset\R^d$ 
\begin{equation}\label{ill605}
	\sup_{t\in [0,T]}\left\| \int_0^t\int_{\!\R^d}\tfrac{\cjap{\xi}^k}{\cjap{\xi_1}^k\cjap{\xi_2}^l}\cos(s|\xi|^2\!-\!s|\xi_1|^2)\cos(s|\xi_2|)\chi_{_{A}}(\xi_1)\chi_{_{B}}(\xi_2)d\xi_1ds\right\|_{_{\!L^2_\xi}} 
	\!\lesssim\ 
	|A|+|B|\,.
\end{equation} 

For $N\in\N$, defining $\vec{N}:=(N,0,\ldots,0)\in \R^d$, 
\begin{equation*}\label{ABR**}
	A_N:=\{\xi_1\in \R^d : |\xi_1|<1/2\}, 
	\quad\quad 
	B_N:=\{\xi_2\in \R^d : |\xi_2-\vec{N}|<1/4\}, 
\end{equation*}
\begin{equation*}\label{tn**}
	R_N:=\{\xi\in \R^d : |\xi-\vec{N}|<1/4\} 
	\quad\quad\textrm{and}\quad\quad 
	t_N:=\frac1{4N^2}\cdot\frac T{1+T}\ , 
\end{equation*} 
then $R_N-B_N\subset A_N$, $\ t_N\in(0,T)$ and, 
for $(\xi_1,\xi_2)\in A_N\times B_N$, 
we have 
\begin{equation*}\label{606} 
\tfrac{\cjap{\xi}^k}{\cjap{\xi_1}^k\cjap{\xi_2}^l}\ \sim\ N^{k-l} 
\quad\quad\textrm{and}\quad\quad 
\cos(s|\xi|^2\!-\!s|\xi_1|^2)\cos(s|\xi_2|)>1/4\,, 
\ \ 
\forall s\in[0,t_N]. 
\end{equation*} 
Thus, from Lemma \ref{lemaleandro*} and \eqref{ill605} yields  
\begin{equation}\label{607} 
t_N|R_N|^{\frac12}|B_N|\,N^{k-l}
\ \lesssim\ 
\left\| N^{k-l}\ \chi_{_{A_N}}\!\!\ast\,\chi_{_{B_N}}\!(\xi)\,\int_0^{t_N}ds\right\|_{_{\!L^2}}
\ \lesssim\ 
|A_N|+|B_N|, 
\quad 
\forall N\in\N.  
\end{equation}
Note that $|A_N|\,,\,|B_N|$ and $|R_N|$ are independent of $N$. 
Hence $l\ \ge\ k-2$ when the \textit{solution mapping} \eqref{datasol} is $C^2$. 
\\ 

	Now we will show that $l\le k+1$. 
	From \eqref{ill01**} follows that \eqref{ill500} holds uniformly for $t\in[0,T]$.
	Let $A,B\subset\R^d$ symmetric sets. 
	By using, in \eqref{ill500}, $\Phi_0=(\varphi,0,0)$ and $\Phi_1=(\upsilon,0,0)$ such that 
	$\varphi,\upsilon\in\mathcal S(\R^d)$, 
	$\langle\cdot\rangle^k\,\widehat{\varphi}\ \sim\ \chi_{_{A}}$ and 
	$\langle\cdot\rangle^k\,\widehat{\upsilon}\ \sim\ \chi_{_{B}}$ we conclude the following estimate 
	for bounded subsets $A,B\subset\R^d$ 
\begin{equation}\label{ill700} 
\sup_{t\in [0,T]}\Biggl\| \int_0^t\!\cos((t-s)|\xi|)|\xi|^2\!\!\int_{_{\R^d}}\!\!\tfrac{\cjap{\xi}^{l-1}}{\cjap{\xi_1}^k\cjap{\xi_2}^k}\cos(|\xi_1|^2s\!-|\xi_2|^2s)\,\chi_{_{_{\!A}}}\!\!(\xi_1)\chi_{_{_{B}}}\!\!(\xi_2)d\xi_1ds   \Biggr\|_{_{\!L^2_\xi}}
\ \lesssim\ 
|A|^{\frac12}|B|^{\frac12}\,. 
\end{equation}

For $N\in\N$, define  
\begin{equation*}\label{ABR***}
	A_N:=\{\xi_1\in \R^d : |\xi_1-\vec{N}|<1/2\}\,\cup\,\{\xi_1\in \R^d : |\xi_1+\vec{N}|<1/2\}, 
	\quad\! 
	B_N:=\{\xi_2\in \R^d : |\xi_2|<1/4\}, 
\end{equation*}
\begin{equation*}\label{tn***}
	R_N:=\{\xi\in \R^d : |\xi-\vec{N}|<1/4\} 
	\quad\quad\textrm{and}\quad\quad 
	t_N:=\frac1{4N^2}\cdot\frac T{1+T}\ . 
\end{equation*} 
Note that $A_N$ and $B_N$ are symmetric. 
Similarly to \eqref{ill605}-\eqref{607}, from \eqref{ill700} we get the following estimate 
\begin{equation*}\label{701} 
t_N|R_N|^{\frac12}|B_N|\,N^{l-k+1}
\ \lesssim\ 
\left\| N^{l-1-k}\ |\xi|^2\ \chi_{_{A_N}}\!\!\ast\,\chi_{_{B_N}}\!(\xi)\,\int_0^{t_N}ds\right\|_{_{\!L^2}}
\ \lesssim\ 
|A_N|^{\frac12}|B_N|^{\frac12}\,, 
\quad 
\forall N\in\N.  
\end{equation*}
Note that $|A_N|\,,\,|B_N|$ and $|R_N|$ are independent of $N$. 
Hence $l\ \le\ k+1$ when the \textit{solution mapping} \eqref{datasol} is $C^2$. 
\hfill $\blacksquare$
%
%
%
%

%
%

\section*{Acknowledgments}

The authors would like to thank Adán Corcho and Xavier Carvajal for their suggestions.

%
%

 
	%
	%

%
%

\end{document}